\documentclass[12pt]{amsart}
\usepackage{amssymb}
\usepackage{euscript}

\let\cal\mathcal

\newtheorem{prop}{Proposition}
\newtheorem{thm}{Theorem}

\newtheorem{lem}{Lemma}
\newtheorem{df}{Definition}

\def\N{\mathbb N}
\def\R{\mathbb R}

\def\dint{\displaystyle\int}

\def\theequation{\thesection.\@arabic\c@equation}

 \makeatletter \@mparswitchfalse \makeatother

\def\eqref#1{(\ref{eq#1})}

\numberwithin{equation}{section}

\begin{document}

\title[THE NEAR RADON-NYKODYM PROPERTY]{THE NEAR RADON-NIKODYM PROPERTY
 IN LEBESGUE-BOCHNER FUNCTION SPACES}
\author{Narcisse Randrianantoanina}
\address{Department of mathematics, University of Texas, Austin, TX 78712}
\curraddr{Department of Mathematics and Statistics, Miami University,
 Oxford, OH 45056}
\email{randrin@muohio.edu}
\author{Elias Saab}
\address{Department of Mathematics, University of Missouri, Columbia, MO
65211}
\email{elias@esaab.math.missouri.edu}

\begin{abstract} Let $X$ be a Banach space and
$(\Omega,\Sigma,\lambda)$  be a finite   measure space, $1 \leq p <\infty$.
It
is shown that $L^p(\lambda,X)$ has the Near Radon-Nikodym property if and
only if $X$ has it. Similarly if $E$ is a K\" othe function space that
does not contain a copy of $c_0$, then $E(X)$ has the Near Radon-Nikodym
property if and only if $X$ does.
\end{abstract}

\maketitle

\section{INTRODUCTION }

 Let X be a Banach space, $(\Omega,\Sigma,\lambda)$  be a finite measure
space, and let $1 \leq p  < \infty$.
We denote by $L^p (\lambda,X) $  the Banach space of all (classes of)
 $\lambda$-measurable  functions from $\Omega$  to $X$
which are $p$-Bochner integrable  with
its usual norm
 $\| f\|_p = \left( \dint \| f(\omega)\|^p \ d\lambda(\omega) \right)^{1/p}$.
 If X is the scalar field then $L^p (\lambda,X) $ will  be
denoted by $L^p (\lambda)$.

The relationship between Radon-Nikodym type properties for Banach spaces and
operators  with domain $L^1[0,1]$ is classical in theory of vector-measures.
Such connections have been investigated by several authors.
In \cite{KPRU}, Kaufman, Petrakis, Riddle and Uhl introduced and studied the
notion of nearly representable operators (see definition below).
 They isolated the class of
Banach spaces $X$ for which every nearly representable operator with range
$X$ is representable. Such Banach spaces are said to have the
Near Radon-Nikodym Property (NRNP). It was shown in \cite{KPRU} that every
Banach lattice that does not contain any copy of $c_0$ has the NRNP; in
particular $L^1$-spaces have the NRNP.
 A question that arises naturally
from this fact is
whether the Lebesgue-Bochner space  $L^1 (\lambda,X) $
   has the NRNP whenever X does.
Let us recall that many related properties such as Radon-Nikodym
property (RNP),
Analytic Radon-Nikodym property (ARNP)
 and complete continuity property (CCP) are known for
Bochner spaces (see \cite{TU}, \cite{DO5} and \cite{RS1} respectively).
We also remark that Hensgen \cite{HE2} observed that (as in the scalar case)
$L^1(\lambda,X)$ has the NRNP if $X$ has the RNP.

  In this paper, we  show that the
Near Radon-Nikodym property can indeed  be lifted from a Banach space $X$ to
 the space $L^1(\lambda,X)$.
 Our proof relies on a representation
of operator from $L^1 $ into $L^1(\lambda,X) $
due to Kalton \cite{KA} and
properties of operator-valued measurable functions along with some well
known characterization of integral and nuclear operators from $L^\infty$
into a given Banach space.

Our  notation is standard Banach space terminology as may be
found in the books \cite{D1}, \cite{DU} and \cite{WO}.

\noindent
{\bf Acknowlegements.} The authors would like to thank
Paula Saab for her constant interests  in this work. The first author
also would like to thank Neal Carothers for creating an enjoyable
work  atmosphere
at the Bowling Green State University  where part of this work was done.

\section{ DEFINITIONS AND PRELIMINARY  RESULTS}

\noindent Throughout this note, $I_{n,k} = [\frac{k-1}{2^n},
\frac{k}{2^n})$ will be the sequence of dyadic intervals in $[0,1]$
 and $\Sigma_n$ is
the $\sigma$-algebra generated by the finite sequence $(I_{n,k})_{k=1,2^n}$.
The word operator will always mean linear bounded operator and $\cal L(E, F)$
will stand for the space of all operators from $E$ into $F$. For any
 given Banach
space $E$, its closed unit ball will be  denoted  by $E_1$.

\begin{df}
 Let X be a Banach space. An  operator
$T: L^1[0,1]\rightarrow X $ is said to be  representable if
 there is a Bochner
integrable function $g \in L^\infty([0,1],X) $
such that T(f)=$\dint fg\
d\mu $ for all f in $L^1 [0,1].$
\end{df}

\begin{df}
  An operator $D: L^1 [0,1]\rightarrow X $
is called a Dunford-Pettis   operator if D sends weakly compact
sets into norm compact sets.
\end{df}

It is well known  (\cite{DU} Example 5,III.2.11) that all representable
operators from
$L^1 [0,1] $ are Dunford-Pettis; but the converse is not true in general.

\begin{df} An operator $T: L^1 [0,1]
\rightarrow X$ is said to be {\bf nearly representable} if for each
Dunford-Pettis operator $D: L^1 [0,1]\rightarrow  L^1 [0,1] $, the
composition   $T\circ D$  is representable.
\end{df}

The notion of nearly representable operators was introduced by Kaufman,
 Petrakis, Riddle and Uhl in \cite{KPRU}.
It should be noted that since the class of Dunford-Pettis operators
from $L^1[0,1]$ into $L^1[0,1]$  is a Banach
lattice (\cite{BO7}),
if an operator $T \in \cal{L}(L^1[0,1],X)$  fails to be nearly
representable then one can find a positive Dunford-Pettis operator $D \in
\cal L(L^1[0,1], L^1[0,1])$ such that $T\circ D$ is not representable.

The following definition isolates the main topic of this paper.
\begin{df}
 A Banach space X has
the {\bf Near Radon-Nikodym Property} ({\bf NRNP}) if every nearly
representable operator from $L^1 [0,1]$ into X is representable.
\end{df}
Examples of Banach spaces with the NRNP are spaces with the RNP,
$L^1$-spaces, $L^1/{H^1}$. For more detailed discussion on the NRNP and
nearly representable operators, we refer to \cite{AP},
\cite{EM6} and  \cite{KPRU}.

We now collect few well known facts about operators from $L^1[0,1]$ that
we will need in the sequel. Our references for these facts are
\cite{BO8}, \cite{BO7} and \cite{DU}.

\noindent {\bf Fact 1.} {\it For a Banach space $X$, there is a one
to one correspondence between  the space of operators from
$L^1[0,1]$ to $X$ and all
 uniformly
bounded
$X$-valued martingales. This correspondence is given by:

\begin{itemize}
 \item [(*)]  $T(f)= \lim\limits_{n \to \infty} \dint \psi_n (t) f(t)\ dt$
if $(\psi_n)_n$ is a uniformly bounded martingale.
 \item [(**)] $\psi_n (t) = 2^n \sum\limits_{k=1}^{2^n} \chi_{I_{n,k}}(t)
\ T(\chi_{I_{n,k}})$ if $T \in \cal L (L^1 [0,1], X)$.
\end{itemize}}

\noindent {\bf Fact 2.} {\it  A uniformly bounded $X$-valued martingale
is Pettis-Cauchy if and only if the corresponding operator $T \in \cal
L(L^1 [0,1],X)$ is Dunford-Pettis.}

As an immediate  consequence of Fact~2, we get:

\noindent{\bf Fact3.} {\it An operator $T \in \cal L (L^1[0,1], X)$
is nearly representable if and only if it maps uniformly bounded
Pettis-Cauchy  martingales to Bochner-Cauchy martingales.}

\begin{df}  Let E and F be Banach spaces and suppose
$T:  E\rightarrow F$ is a bounded linear operator. The operator $T$ is
said to be an
{\bf absolutely summing operator } if there is a constant C
such that for
any finite sequence $(x_m)_{1 \leq  m \leq n} $ in E, the following
holds:
$$\sum_{m=1}^n ||T x_m || \leq C \  \sup \left \{ \sum_{m=1}^n
| x^* (x_m)|;\    x^*\in E^* \ ; ||x^*||\leq 1 \right\}.$$
\end{df}

The least constant $C$ for the inequality above to
hold will be denoted by $\pi_1( T)$.
It is well known  that
the class of all absolutely summing operators from E to F is a Banach
space under the norm $\pi_1(T)$. This Banach space will
be denoted by $\Pi_1 (E,F)$.

\begin{df} We say that an operator $T: E \rightarrow F$ is an
{\bf integral operator} if it admits a factorization:
$$\begin{array}{ccc}
 E& \stackrel{i\circ T}{\longrightarrow}&  F^{**} \cr
\Big\downarrow\vcenter{%
\rlap{$\alpha$}}&
&\Big\uparrow\vcenter{%
\rlap{$\beta$ }} \cr
L^\infty(\mu)& \stackrel{J}{\longrightarrow}&  L^1(\mu)
\end{array}$$
where $i$ is the inclusion from $F$ into $F^{**}$, $\mu$ is a probability
measure on a compact space
$K$,
$J$ is the natural inclusion and
$\alpha$ and $\beta$ are bounded linear operators.
\end{df}

We define the integral norm $i(T)= \inf\{||\alpha||. ||\beta||\}$
where the infimum is taken over all such factorization.
We denote  by $I(E,F)$ the space of integral operators from $E$ into $F$.

If E = C(K) where K is a compact Hausdorff
space or $E=L^\infty(\mu)$ then it is well known that T is absolutely
summing (equivalently $T$ is integral) if and only if its
representing measure G (see \cite{DU}, p.152) is of bounded variation
and in this case  $\pi_1(T)= i(T) = | G|(K)$
where $ |G|(K)$ denotes the total variation of $G$.

\begin{df} We say that an operator $T: E \rightarrow  F$ is a {\bf
nuclear operator} if there exist sequences $(e_{n}^*)_n$
in $E^*$ and $(f_n)_n$ in $F$ such that $\sum\limits_{n=1}^\infty
||e_{n}^*||\ ||f_n|| < \infty$ and such that
   $$T(e)=\sum_{n=1}^\infty e_{n}^*(e)f_n $$
for all $e \in E$.
\end{df}
We define the nuclear norm
$n(T)=\inf\{\sum\limits_{n=1}^\infty ||e_{n}^*||\ ||f_n|| \}$ where the
infimum is taken over all sequences $(e_{n}^*)_n$ and $(f_n)_n$ such that
$T(e) = \sum\limits_{n=1}^\infty e_{n}^*(e)f_n$ for all $e \in E$. We
denote by $N(E,F)$ the space of all nuclear operators from $E$ into $F$
under the norm $n(.)$.

\noindent
{\bf Fact 4.} {\it An operator $T \in \cal L(L^1[0,1],X)$ is
representable
if and only if its restriction to $L^\infty[0,1]$, $T|_{L^\infty[0,1]}
\in\cal L(L^\infty[0,1],X)$ is nuclear.}

Throughout this paper, we will identify the two function spaces
$L^p(\lambda,L^p(\mu,X))$ and  $L^p(\lambda \otimes \mu,X)$ for $1 \leq
p < \infty$ (see \cite{DS}, p.198).

The following representation theorem of Kalton \cite{KA}  is
essential for the proof of the  main result.
We denote by $\beta (K)$   the $\sigma  $-Algebra of
Borel subsets of K in the statement of the theorem.

\begin{thm}\cite{KA}(Kalton)
  Suppose that:
\begin{itemize}
 \item  [(i)] K is a compact metric space and  $ \mu $ is a Radon
probability measure on K;
 \item  [(ii)] {\bf $\Omega $ } is a Polish space and $\lambda $ is a
Radon measure on
 $\Omega $;
 \item  [(iii)] $ X$  is a separable Banach space;
\item   [ (iv)] $T:  L^1 (\mu) \longrightarrow  L^1 (\lambda ,X ) $ is a
bounded linear operator.
\end{itemize}
  Then there is a map  $ \omega \to   T_\omega  \   (\Omega
\rightarrow \Pi_1(C(K),X))$ such that for every f $\in$
C(K),
the map $ \omega   \to  T_ \omega (f) $ is Borel measurable from
$\Omega $     to X and
\begin{itemize}
\item [($\alpha $)]  If $\mu_ \omega$ is the representing measure of $ T_
\omega $
then  $$\int_{\Omega } |\mu_\omega|(B)\ d\lambda(\omega)
\leq  ||T||\mu (B) \quad \text{for every}\ B \in\beta (K);$$
\item [($\beta $)]  If $f\in L^1 (\mu ) $,  then for $\lambda$ a.e
$\omega $, one has $ f \in L^1 (|\mu_\omega |)$;
\item [($ \gamma $)] $Tf(\omega) =  T_\omega (f) $    for $\lambda$ a.e
$\omega$ and for every $ f \in  L^1 (\mu )$.
\end{itemize}
\end{thm}
  The following proposition gives a characterization of representable
operators in connection with  Theorem~1.

\begin{prop}(\cite{RS3})
  Under the assumptions of Theorem 1, the
following two statements are equivalent:
\begin{itemize}
\item [(i)] The operator $T$ is representable;
\item [(ii)] For\  $\lambda$ a.e $\omega$, $\mu_\omega$
 has Bochner integrable density with respect to $\mu $.
\end{itemize}
\end{prop}

\noindent For the next result, we need the following definition.

\begin{df}
 Let $E$ and $F$ be Banach spaces.
 A map $T:(\Omega, \Sigma, \lambda)
\rightarrow
\cal L(E,F)$ is said to be strongly measurable if $\omega
\rightarrow T(\omega)e$ is measurable for every $e \in E$.
\end{df}

 We observe that  if $E$ and $F$ are separable
Banach spaces and $T:(\Omega,\lambda)\rightarrow \cal L(E,F)$ with
$\sup\limits_{\omega} \Vert T(\omega) \Vert \leq 1 $, then $T$ is
strongly measurable if and only if $T^{-1}(B)$ is $\lambda$-measurable
for each Borel subset $B$ of $\cal L(E,F)_1$ endowed with the stong
operator topology.

\noindent
The following selection result will be needed for the proof of the main
theorem.

\begin{prop}
 Let $X$ be a separable Banach space and
 $T: (\Omega,\lambda) \rightarrow \cal L (L^1 [0,1],X)$ be a strongly
measurable map with:
\begin{itemize}
 \item [(1)]  $\sup\limits_{\omega} \Vert T(\omega) \Vert \leq 1$;
 \item [(2)]  $T(\omega)$ is not nearly representable for $\omega \in A$,
$\lambda(A) >0$.
\end{itemize}
Then one can choose a map $D: (\Omega,\lambda) \rightarrow \cal L(L^1
[0,1], L^1[0,1])$ strongly measurable such that:
\begin{itemize}
 \item [(i)] $\sup\limits_{\omega} \Vert D(\omega) \Vert \leq 1 $;
 \item [(ii)] $ T(\omega) \circ D(\omega)$ is not representable for each
$\omega \in A$;
 \item [(iii)] $D(\omega)$ is  Dunford-Pettis for a.e $\omega \in \Omega$;
\item [(iv)] $D(\omega)$ is a positive operator for every $\omega \in
\Omega$.
\end{itemize}
\end{prop}

We will need several steps for the proof.

\begin{lem} The space $\cal L(L^1[0,1],X)_1$, the unit
ball of the space $\cal L(L^1[0,1],X)$ endowed with the strong operator
topology is a Polish space.
\end{lem}

\begin{proof} Let us consider the Polish space $\Pi_n \{
X^{2^n}\}$. We will show that $\cal L (L^1[0,1],X)_1$ is homeomorphic to
a closed subspace of $\Pi_n \{X^{2^n}\}$.

Let $\cal C$ be the following subset of $\Pi_n \{X^{2^n}\}$:
$(x_{n,k})_{k \leq 2^n ; n \in \N}$ belongs to $\cal C$ if and only if
\begin{itemize}
   \item [(a)] $x_{n,k} = \frac{1}{2} ( x_{n+1, 2k-1} + x_{n+1,2k})$ for
all $k \leq 2^n $ and $n \in \N$;
   \item [(b)] $\Vert x_{n,k} \Vert \leq 1$ for all $k \leq 2^n$ and $n
\in \N$.
\end{itemize}

It is evident that $\cal C$ is closed in $\Pi_n
\{X^{2^n}\}$.

\noindent
Consider the map $\Gamma: \cal L (L^1[0,1],X)_1 \rightarrow \Pi_n
\{X^{2^n}\}$ given by
     $T \rightarrow  (2^n T(\chi_{I_{n,k}}))_{k \leq 2^n, n \in
\N}.$

\noindent The map
$\Gamma$ is clearly continuous, one to one
 and its range is contained in $\cal C$.
We claim that $\Gamma(\cal L(L^1[0,1],X)_1)= \cal C$ and $\Gamma
\vert_{\cal C}^{-1}$ is continuous:
to see this claim,
 let  $x=(x_{n,k}) \in \cal C$ and  $T\in \cal L(L^1[0,1], X)$
defined by the martingale $ \psi_n(t) = \sum\limits_{k=1}^{2^n} x_{n,k}
\chi_{I_{n,k}}(t)$. The operator $T$ is well defined  (see Fact 1) and
$T(\chi_{I_{n,k}})=
(1/2^n)  x_{n,k}$ so $\Gamma(T)=x$. Using the fact that span $\{
\chi_{I_{n,k}}, k \leq 2^n,
n \in \N\}$ is dense in $L^1[0,1]$, the continuity of $\Gamma\vert_{\cal
C}^{-1}$ follows. The lemma is proved.
\end{proof}

 Consider $\cal L(L^1[0,1],X)_1$ with the strong operator
topology and $L^1([0,1],L^1[0,1]) $ with   the norm-topology.

Using the fact that  the natural injection from
$L^\infty([0,1],L^1[0,1])$ into
$L^1([0,1],L^1[0,1])$ is a semi-embedding, the unit ball of
$L^\infty([0,1],L^1[0,1])$ (that
we will denote by $Z$) is a closed subset of the
Polish space  $L^1([0,1],L^1[0,1])$ so $Z$
with the relative topology is a Polish space.

The space $\cal L(L^1[0,1],X)_1 \times Z^\N$ with the product topology
is a Polish space.

\noindent Let $\cal A $ be a subset of $\cal L(L^1[0,1],X)_1 \times
Z^\N$   defined  as follows:

$\{T, (\phi_n)_n \} \in \cal A $ if and only if
\begin{itemize}
 \item [(i)] $\mathbb{E}(\phi_{n+1}/ \Sigma_n) = \phi_n $
 for every $n \in \N$;
\item [(ii)] $\lim\limits_{n,m} \sup\limits_{g \in L^\infty , \Vert g
\Vert_\infty \leq 1} \dint \vert \dint (\phi_m (t,s)- \phi_n(t,s))g(s) \
ds
\vert \ dt =0$;
\item [(iii)] $\lim\limits_{j \to \infty} \sup\limits_{n,m \geq j} \dint
\Vert T(\phi_n(t)
-\phi_m(t))
\Vert \ dt >0$;
\item [(iv)] $\phi_n \geq 0$ as element of the Banach lattice
$L^\infty([0,1], L^1[0,1])$.
\end{itemize}

\begin{lem} The set $\cal A$ is a Borel subset of $\cal L
(L^1[0,1],X)_1 \times Z^\N$.
\end{lem}

\begin{proof} $(i)$ Let $\cal A_1$ be a subset of
$Z^\N$ given by  $\phi =(\phi_n)_n \in \cal A_1 $ if
and only if
 $$\mathbb{E}(\phi_{n+1}/ \Sigma_n )=\phi_n \quad \forall n \in \N. $$
We claim that
 $\cal A_1$ is a Borel subset of $Z^\N$: if we denote by $P_n$ the
$n^{th}$ projection of $Z^\N$ and $\mathbb{E}_n$ the conditional expectation
 with respect to
$\Sigma_n$, then the map
  $\theta_n: L^1([0,1],L^1[0,1])^\N \rightarrow
L^1([0,1],L^1[0,1])$
given by $\theta_n (\phi) =(\mathbb{E}_n \circ P_{n+1} -P_n)(\phi)$
 is continuous
and therefore
$\cal A_1 = \bigcap\limits_{n\in \N} \theta_n^{-1} (\{0\}) \cap Z^\N$
is Borel measurable.

\noindent $(ii)$  Let $g \in L^\infty $ be fixed. For  every $m ,n \in \N$,
the map:
\begin{equation*}
\begin{split}
   L^1([0,1],L^1[0,1])^\N  &\longrightarrow \R \cr
     \phi   &\longrightarrow \int \vert \int (\phi_m (t,s)-\phi_n(t,s))
g(s)\ ds \vert \ dt
\end{split}
\end{equation*}

is continuous so
$ \phi \rightarrow \Gamma_{n,m} (\phi) =\sup_{g \in L^\infty , \Vert
g \Vert
\leq 1} \dint \vert \dint \phi_m (t,s)- \phi_n(t,s))g(s) \ ds \vert \ dt $
is lower semi-continuous and therefore
 $\phi \rightarrow \Gamma (\phi)= \lim_{j \to \infty}\sup_{n,m \geq
j}
\Gamma_{n,m}(\phi)$ is Borel measurable and  we have that
$$\cal A_2 = \{ \phi:\ \lim_{n,m} \sup_{g \in L^\infty, \Vert g \Vert
\leq 1} \int \vert \int (\phi_m(t,s)- \phi_n(t,s) )g(s)\ ds \vert \ dt
 =0\} \cap Z^\N $$ is a Borel measurable subset of
$Z^\N$.

\noindent $(iii)$ For each $n$ and $m$ in $\N$, the map
\begin{equation*}
\begin{split}
  \theta_{n,m}: \cal L(L^1[0,1],X)_1 \times
L^1([0,1],L^1[0,1])^\N &\longrightarrow \R \cr
(T, \phi) &\longrightarrow \int \Vert T(\phi_n(t)) -T(\phi_m(t)) \Vert \ dt
\end{split}
 \end{equation*}
is continuous and then the set
$\cal B = \{(T,\phi)   ; \ \limsup\limits_{n,m} \theta_{n,m} (T,\phi) >0
\}$ is a Borel measurable subset of  $\cal L(L^1[0,1],X)_1 \times
L^1([0,1],L^1[0,1])^\N$.

\noindent $(iv)$ The set $\cal P$ of sequences of positive functions is a
closed subspace of $Z^\N$ .

Now $\cal A = \cal B \cap \{ \cal L(L^1[0,1],X)_1 \times (\cal A_1 \cap
\cal A_2 \cap \cal P) \} $ so $\cal A$ is Borel measurable. The lemma is
proved.
\end{proof}

\noindent
{\it Proof of Proposition~2.} Let  $U$  be the
restriction on $\cal A$ of the first projection; the set
$U(\cal A)$ is an analytic subset of $\cal L(L^1[0,1],X)_1$ and by
Theorem 8.5.3 of \cite{CO}, there is an universally measurable map
$\theta: U(\cal A) \rightarrow Z^\N $
such that  the graph of $\theta$ is
contained in  $\cal A$.

By assumption, $T: (\Omega, \lambda) \rightarrow \cal L(L^1([0,1],X)_1$
 is Lusin-measurable for the strong operator topology
and $T(\omega) \in U(\cal A)$ for evry  $\omega \in A$. So the following map
   \begin{equation*}
   \begin{split}
       \Omega  &\longrightarrow L^1([0,1],L^1[0,1])^\N \cr
       \omega  &\longrightarrow  \begin{cases} \theta(T(\omega))   \quad
&\text{if $\omega \in A$}\cr
0 \quad &\text{otherwise}
\end{cases}
\end{split}
 \end{equation*}
is well-defined and is $\lambda$-measurable. Moreover for every
 $\omega \in A$,
 $\{ T(\omega), \theta(T(\omega)) \} \in \cal A. $

Let $Q_n$ be the $n^{\text{th}}$ projection from $Z^\N$ onto $Z$. For
every $n \geq 1 $, let $\phi_n (\omega) = Q_n ( \theta(T(\omega))$. By
construction, the sequence $(\phi_n(\omega))_n$ is a uniformly bounded
martingale from $[0,1]$ into $L^1[0,1]$, so it defines an operator from
$L^1[0,1]$ into $L^1[0,1]$ by
   $$D(\omega)(f) = \lim_{n \to \infty} \dint \phi_n (\omega)(t)f(t)\
dt.$$
Notice that for every $f \in L^1[0,1]$,  the map
 $M_f: Z \rightarrow
L^1([0,1],L^1[0,1])$ defined by $M_f (h) =f.h$ is continuous and
$D(\omega)(f) = \lim\limits_{n \to \infty} \dint M_f(Q_n (\theta(T(\omega)))\
dt$.
We conclude that for every $f \in L^1[0,1]$, the map
         $ \omega \rightarrow D(\omega)(f)$ ($\Omega \to L^1[0,1]$)
is measurable.  Now
condition $(iii)$ implies that $T(\omega) \circ D(\omega)$ is not
representable for $\omega \in A$
and condition $(iv)$ insures that $D(\omega) \geq 0$
 for every $\omega \in \Omega$.
\qed

The following proposition is crusual for the proof our main result and could
be of independent interest.

\begin{prop}
Let $\omega \to D(\omega)$  $(\Omega \to \cal L(L^1[0,1], L^1[0,1])_1)$ be a
strongly measurable map such that $D(\omega)$ is positive
and Dunford-Pettis for every
$\omega \in \Omega$. If we denote by $\theta(\omega)$ the restriction of
$D(\omega)$ on $L^\infty[0,1]$, then $\omega \to \theta(\omega)$ is
norm-measurable as a map from $\Omega$ into $I(L^\infty[0,1],L^1[0,1])$.
\end{prop}

We will begin by proving  the following simple lemma.

\begin{lem}
Let $D: L^1[0,1] \to L^1[0,1]$ be a positive Dunford-Pettis operator and
$\theta=D|_{L^\infty}$. Then $\theta$ is compact integral and is weak$^*$
to weakly  continuous. Moreover $i(\theta)=||\theta||$.
\end{lem}

\begin{proof} The fact that $\theta$ is compact
integral is trivial. For the weak$^*$ to weak continuity, we observe that
$\theta^*(L^\infty[0,1]) \subset L^1[0,1]$.
 For the identity of the norms, we will use the fact that
$i(\theta)$ is equal to the total variation of the representing measure
of $\theta$.

Let $G$ be the representing measure of $\theta$ and $\pi$ be a finite
measurable partition of $[0,1]$. We have the following:
\begin{equation*}\begin{split}
\sum_{A \in \pi} ||G(A)||_{L^1} &= \sum_{A \in \pi} ||D(\chi_A)|| \cr
&\leq \sum_{A\in \pi} ||\ |D|(\chi_A)\ || \cr
&= \sum_{A \in \pi} ||\ |\theta|(\chi_A)\ || \cr
&= \sum_{A \in \pi} \int |\theta|(\chi_A)(t)\ dt \cr
&=\int |\theta|(\chi_{[0,1]})(t)\ dt  \leq ||\ |\theta|\ ||
\end{split}\end{equation*}
where $|D|$ and $|\theta|$ denote the modulus of $D$ and $\theta$
respectively (see \cite{LT}).
So by taking the supremum over all finite measurable partition of [0,1], we get
that $i(\theta) \leq ||\ |\theta|\ ||$ and since $\theta$ is a positive
operator, $|\theta|=\theta$. The lemma is proved.
\end{proof}
\noindent
{\it Proof of  Proposition~3.} Notice that $\theta(\omega) \in
K_{w^*}(L^\infty[0,1],L^1[0,1])$ for every $\omega \in \Omega$ where
$K_{w^*}(L^\infty[0,1],L^1[0,1])$ denotes the space of compact operators
from $L^\infty[0,1]$ into $L^1[0,1]$ that are weak$^*$ to weakly
continuous. So we get that $\omega \to \theta(\omega)$ is strongly
measurable and is separably valued ($K_{w^*}(L^\infty[0,1],L^1[0,1]) =
L^1[0,1] \widehat{\otimes}_\epsilon L^1[0,1]$ where
$\widehat{\otimes}_\epsilon$ is the injective tensor product). By the
Pettis measurability theorem (see Theorem II-2 of \cite{DU}), the map
$\omega \to \theta(\omega)$ is measurable
for the norm operator topology.

For each $n \in \N$, let  $\mathbb{E}_n$ be the conditional expectation operator
with respect to $\Sigma_n$. The sequence $(\mathbb{E}_n)_n$
 satisfies  the following
properties:
 $(\mathbb{E}_n)_n$  is a
sequence of finite  rank operators
in $\cal L(L^1[0,1], L^1[0,1])_1$,
 $\mathbb{E}_n \geq 0$
for every $n \in \N$ and $(\mathbb{E}_n)_n$ converges to
 the identity operator $I$ for
the strong operator topology. Consider $S_n= \mathbb{E}_n \wedge I$.
 Since $S_n \leq \mathbb{E}_n$ and $\mathbb{E}_n$ is
 integral (it is of finite rank), one can deduce from
 Grothendieck's characterization of integral operators
with values in  $L^1[0,1]$ (see
for instance \cite{DU}  p. 258) that $S_n$ is also integral.

\noindent{\bf Sublemma.} {\it For each $n \in \N$,
 there exists $K_n \in \text{conv}\{S_n, S_{n+1},
\dots\}$ such that the sequence $(K_n)_n$ converges to $I$ for the strong
operator topology.}

For this, we observe first that $(S_n(f))_n$ converges weakly to $f$ for every
$f \in L^1[0,1]$; in fact, if $f\geq 0$  and $n \in \N$ then
$S_n(f)=\inf\{\mathbb{E}_n(g) + (f-g);\ 0\leq g \leq f\}$. Choose
 $0 \leq g_n \leq  f$ such that
 $\|S_n(f)-(\mathbb{E}_n(g_n) +(f-g_n))\|_1 \leq 1/n$. Since $[0,f]$ is weakly
compact, we can assume (by taking a subsequence if necessary) that $(g_n)_n$
converges weakly to a function $g$. To conclude that $S_n(f)$ converges
weakly, notice that if $\varphi \in L^\infty[0,1]$ then
 $\lim\limits_{n \to \infty}\mathbb{E}_{n}^*(\varphi)= \varphi$ a.e
 ($\mathbb{E}_{n}^*=\mathbb{E}_n$). So we have for every $n \in \N$,
$|\langle S_n(f)-f,\varphi\rangle| \leq 1/n +
 |\langle\mathbb{E}_n(g_n) - g_n, \varphi\rangle|$ and
$$
 |\langle \mathbb{E}_n(g_n) - g_n, \varphi\rangle| =
 |\langle  g_n, \mathbb{E}_n(\varphi) -\varphi\rangle|
 \leq\langle  f , | \mathbb{E}_n(\varphi) -\varphi |\rangle.
$$
And by the Lebesgue dominated convergence, we have
$\lim\limits_{n \to \infty}\langle\mathbb{E}_n(g_n)- g_n, \varphi\rangle=0$.
Now fix $(f_k)_k$ a countable dense subset of  the closed unit ball of
$L^1[0,1]$.  For $k=1$, we can choose by Mazur's theorem,  a sequence
$(S_{n}^{(1)})_n$ with $S_{n}^{(1)} \in \text{conv}\{ S_n, S_{n+1}, \dots\}$
for every $n \in \N$ and such that
 $\lim\limits_{n \to \infty}|| S_{n}^{(1)}(f_1)- f_1 ||=0$. By induction,
one can use the same argument to construct $S_{n}^{(k+1)} \in
\text{conv}\{ S_{n}^{(k)}, S_{n+1}^{(k)}, \dots\}$ such that
$\lim\limits_{n\to \infty} || S_{n}^{(k+1)}(f_j) -f_j ||=0$ for every
 $j \leq (k+1)$. From Lemma~1 of \cite{T2}, one can fix a sequence
$(K_n)_n$ such that for every $k \in \N$, there exists  $n_k \in \N$ such that
for $n \geq n_k$, $K_n \in \text{conv}\{S_{n}^{(k)}, S_{n+1}^{(k)}, \dots\}$.
From this, it is clear that $\lim\limits_{n \to \infty}|| K_n (f_k) -f_k ||=0$
for every $k \in \N$ and since $(f_k)_k$ is dense and
 $\sup\limits_{n}|| K_n||\leq 1$, $(K_n)_n$ verifies the requirements of
the sublemma.

To complete the proof of the proposition, let $(K_n)_n$
 be as in the above sublemma and consider
  $C_n: K_{w^*}(L^\infty[0,1], L^1[0,1]) \rightarrow I(L^\infty[0,1],
L^1[0,1])$ ($T \to K_n \circ T$). Since $K_n$ is integral, the map
$C_n$ is well-defined and is clearly continuous. Therefore
$\omega \to K_n \circ \theta(\omega)$ is measurable for the integral norm.
Since $(K_n)$ converges to $I$ for the strong operator topology and
$\theta(\omega)$ is compact,
 $\lim\limits_{n \to \infty}\| K_n \circ \theta(\omega) - \theta(\omega)\|=0$.
Observe that $K_n\circ \theta(\omega)\leq \theta(\omega)$ for
every  $\omega \in\Omega$ and  for every  $n \in \N$.
 We conclude from Lemma~3 that
 $i(\theta(\omega)- K_n\circ \theta(\omega))=\| \theta(\omega)- K_n \circ
\theta(\omega)\|$ and hence for a.e $\omega \in \Omega$,
$$\lim_{n \to \infty } i(\theta(\omega ) - K_n \circ \theta(\omega))=0.$$
and since $K_n\circ \theta(.)$'s are  measurable
so is  $\theta(.)$, the proposition is
proved.
\qed

The following proposition is probably known but we do not know of any
specific reference.

\begin{prop} Let $X$ be a Banach space
and $S:(\Omega,\lambda)
\rightarrow \cal L(L^1[0,1],X)$ be a strongly measurable map with
$\sup\limits_{\omega} \Vert S(\omega) \Vert \leq 1$. Then the following
assertions are equivalent:
\begin{itemize}
 \item [(a)] The operator
    $ H: L^1(\Omega\times [0,1],\lambda \otimes m)
\rightarrow X$
given by
 $ H(f )= \dint_{\Omega} S(\omega)(f(\omega ,.))\
d\lambda(\omega)$
 is representable;
 \item [(b)] The operator
$ K: L^1[0,1] \rightarrow L^1(\lambda,X)$ given by
          $  K(g)= S(.)g $ is representable;
 \item [(c)] $S(\omega)$ is representable for a.e $\omega \in \Omega$.
 \end{itemize}
 \end{prop}

\begin{proof} $(a) \Rightarrow (b)$\  If $H$ is representable, then we
can find an essentially bounded measurable map $\psi: \Omega \times [0,1]
\rightarrow X$ that represents $H$. The map $\psi^\prime: [0,1]
\rightarrow L^1(\lambda,X)$ given by $t \rightarrow \psi(.,t)$
belongs to $L^\infty([0,1],L^1(\lambda,X))$: in fact
$||\psi^\prime(t)||=\dint_{\Omega} ||\psi(\omega,t)||\ d\lambda(\omega)$
for every $t \in [0,1]$ hence $||\psi^\prime||_\infty \leq
||\psi||_\infty$ and we claim  that $\psi^\prime$ represents
$K$: for each
$g \in L^1[0,1]$,
$\{\dint
\psi^\prime (t)g(t)\ dt\}(\omega)=\dint \psi(\omega,t)g(t)\ dt$ for
a.e $\omega$. For every measurable subset $A$ of $\Omega$,
\begin{equation*}\begin{split}
   \int_A Kg(\omega)\ d\lambda(\omega) &= H(\chi_A \otimes g) \cr
    &= \iint \psi(\omega,t)g(t) \chi_A (\omega)\ dt\ d\lambda(\omega) \cr
    &= \int_A \{\int \psi^\prime (t) g(t)\ dt \}(\omega)\
d\lambda(\omega)
\end{split}\end{equation*}
which shows that $Kg =\dint \psi^\prime (t) g(t)\ dt$.

 $(b) \Leftrightarrow (c)$ Let
$\mu_\omega \in M([0,1],X)$ be
the representing measure for $S(\omega)$ (i.e
 $S(\omega)(\chi_A) = \mu_\omega(A)$). It is well known  that
$S(\omega)$ is representable if and only if $\mu_\omega $ has Bochner
density with respect to $dt$. Notice now that $K(g) (\omega)
=S(\omega)(g)=
 \dint g(t)\ d\mu_\omega(t)$. Hence, by the uniqueness of the
representation of Theorem 1 (see \cite{KA}, p.316), the
family
$(\mu_\omega)_\omega$ represents $K$.
Apply now Propostion~1 to conclude the equivalence.

\noindent $(b) \Rightarrow (a)$ If $\psi^\prime: [0,1] \rightarrow
L^1(\lambda,X)$ represents $K$, then there is a
 map $\Gamma:\Omega\times[0,1]\rightarrow X$ so that $\Gamma \in
L^1(\Omega \times [0,1],\lambda \otimes m)$
      and      $\Gamma(.,t)=
\psi^\prime(t)$ for a.e
$t\in [0,1]$
(see \cite{DS}, p.198).
 We
claim that $\Gamma \in L^\infty (\Omega \times [0,1],\lambda \otimes
m)$ and
 represents
$H$. To prove this claim, let $G(V)= H(\chi_V)$ be the representing
measure of
$H$. If
$A$ is
 a measurable subset of
$\Omega$ and $I$ is a measurable subset of $[0,1]$, we have the following:
 \begin{equation*}\begin{split}
   G(A \times I) &= H(\chi_A \otimes \chi_I ) \cr
                  &=\int_\Omega K(\chi_I) \chi_A\
d\lambda(\omega)\cr
        &=\int_A (\int_I \psi^\prime(t)\ dm(t))(\omega)\
d\lambda(\omega)\cr
        &=\iint_{A \times I} \Gamma (\omega,t) \ d(\lambda \otimes m)
(\omega,t).
        \end{split}\end{equation*}
This will imply that $G(V) = \displaystyle\iint_V
\Gamma (\omega,t) \ d(\lambda \otimes m)
(\omega,t) $ for every Borel subset of $\Omega \times [0,1]$. Apply now
Lemma 4-III of \cite{DU} to conclude that $H$ is representable.
\end{proof}

\section{MAIN RESULT }

\begin{thm}
 Let X be a Banach space and $(\Omega ,\Sigma ,
\lambda)$
 be  a finite measure space then  $L^1 (\lambda ,X ) $ has the NRNP
if and only if  X does.
\end{thm}

 For the proof, let us assume without loss of generality that X is
seperable, $\Omega$ is a compact metric space and $\lambda $ is a Radon
measure in the Borel $\sigma$-Algebra  $\Sigma $ of $\Omega $. For what
follows, $J_X$ denotes the natural inclusion from $L^\infty(\lambda,X$ into
$L^1(\lambda,X)$.

We will begin with  the proof of the following special case.
\begin{prop}
Let $X$ be a Banach space with the NRNP and
$T: L^1[0,1] \to L^\infty(\lambda,X) $ be a bounded linear operator.
Then $J_X\circ T$ is representable if and only if it is nearly
representable.
\end{prop}
\begin{proof}
Let $T: L^1[0,1] \rightarrow L^\infty(\lambda,X)$ be a bounded
operator with $||T|| \leq 1$.
 By Lemma~1 of \cite{RS1}, there exists a strongly measurable map
 $\omega \to T(\omega)\ (\Omega \to \cal L(L^1[0,1],X)_1)$ such that
 $Tf(.)= T(.)f$ for every $f \in L^1[0,1]$.

 Assume that $J_X\circ T$ is nearly representable but not representable.
 Proposition~4  asserts that there exists a measurable subset $A$ of
$\Omega$
 with $\lambda(A)>0$ and such that $T(\omega)$ is not representable for each
 $\omega \in A$. Since $X$ has the NRNP, the operator $T(\omega)$ is not
nearly representable for
 each $\omega \in A$. Using our selection result (Proposition~2), one
 can choose a strongly measurable map $\omega \to D(\omega) \ (\Omega \to
\cal L(L^1[0,1], L^1[0,1])_1)$
 such that $D(\omega)$ is positive,  Dunford-Pettis for every $\omega
\in
\Omega$
 and $T(\omega)\circ D(\omega)$ is not representable for every
 $\omega \in A$. It should be noted that if
  $D \in \cal L(L^1[0,1], L^1[0,1])$ is  a Dunford-Pettis operator,
 since $J_X\circ T$ is nearly representable,  we get that  $T(\omega) \circ
D$ is
 representable for a.e $\omega \in \Omega$ (see Proposition~4). However
the exceptional
  set may depend  on the operator $D$.

  As before let $\theta(\omega)= D(\omega)|_{L^\infty}$.
We deduce from  Proposition~3 that
  the map
 $\omega \to \theta(\omega)\ (\Omega \to I(L^\infty[0,1],L^1[0,1]))$ is
norm-measurable.

Let $(\Pi_n)_{n \in \N}$ be a sequence of finite measurable partition
of $\Omega$ such that $\Pi_{n+1}$ is finer than $\Pi_n$ for every $n \in \N$
and $\Sigma$ is generated by $\bigcup\limits_{n\in \N}\{B\ ; B\in \Pi_n \}$.

For each $B \in \Sigma$, we denote by $D_B$ the operator defined as follows:

\noindent
$D_B (f)= \dint_B D(\omega)(f)\ d\lambda(\omega)$ and define
$D_n(\omega) = \sum_{B \in \Pi_n} \frac{D_B}{\lambda(B)}\ \chi_{B}(\omega)$.
The operator $D_B$ is a Dunford-Pettis operator for each $B
\in \Sigma$
(see \cite{VO} Theorem 1.3) and therefore
$D_n(\omega)$ is Dunford-Pettis for each
$n\in\N$ and $\omega \in \Omega$.

\noindent Claim: The operator $T(\omega)\circ D_n(\omega)$ is representable
for a.e $\omega \in \Omega$.

To see this  claim, notice that $T(\omega)\circ D_B$ is representable
for a.e $\omega \in \Omega$. Fix a set $N_B$ with $\lambda(N_B)=0$
such that $T(\omega)\circ D_B$ is representable for $\omega \notin N_B$;
let $ N = \bigcup\limits_{n \in \N} \bigcup\limits_{B \in \Pi_n} N_B$;
 $\lambda(N)=0$ and for $\omega \notin N$, we have
 $T(\omega)\circ D_n(\omega) = \sum\limits_{B \in \Pi_n} \frac{T(\omega)
 \circ D_B}{\lambda(B)}\ \chi_B(\omega)$ is representable.

 Now if we denote by $\theta_n$ (resp. $\theta_B$) the restriction
 on $L^\infty [0,1]$ of $D_n$ (resp. $D_B$), we have
 $$ \theta_n(\omega)= \sum_{B \in \Pi_n} \frac{\theta_B}{\lambda(B)}
 \chi_B(\omega)$$
 for each $\omega \in \Omega$, and since $\theta(.)$ is norm-measurable
(see Proposition~3), we get that
 $$ \theta_n(\omega)= \sum_{B \in \Pi_n} \frac{\text{Bochner}-\dint_B
 \theta(s)\ d\lambda(s)}{\lambda(B)} \ \chi_B (\omega).$$
 It is well known (see for instance \cite{DU} Corollary V-2 ) that
 $\theta_n(.)$ converges (for the integral norm) to $\theta(.)$ a.e.
 Now since $T(\omega)\circ D_n(\omega)$ is representable for a.e
$\omega$,
 the operator $T(\omega)\circ \theta_n(\omega)$ is nuclear for a.e
 $\omega$ and since $\theta_n(\omega)$ converges a.e to
 $\theta(\omega)$ for the integral norm, we get that
 $$ \lim_{n \to \infty} i\left( T(\omega)\circ \theta_n(\omega) -
 T(\omega)\circ \theta(\omega) \right)=0 \quad \text{for a.e} \
 \omega \in \Omega.$$
 As a result, the operator $T(\omega)\circ \theta(\omega)$ is nuclear
 for a.e $\omega \in \Omega$ and this is equivalent to
 that $T(\omega)\circ
 D(\omega)$ being representable for a.e $\omega \in \Omega$. Contradiction.
\end{proof}

For the general case,
let $T:  L^1 [0,1] \rightarrow  L^1(\lambda ,X ) $ be  a nearly
representable operator and fix a strongly Borel measurable  map
\  $\omega \rightarrow T_\omega$   ($\Omega \rightarrow
\Pi_1(C[0,1],X)$)
 as in Theorem~1.   Let us denote by $\mu_\omega $ the representing
measure of $T_\omega $.
Our goal is to show that for $\lambda $ a.e
$\omega  , \,\mu_\omega$ has a Bochner integrable density with respect
to the Lebesgue measure $m$ in $[0,1]$. This will imply that $T$ is
representable by Proposition 1. To do that, we need to establish
several steps:

\begin{lem}
For $\lambda \ \text {a.e} \quad  \omega $  in $\Omega $, we have
 $\vert \mu_\omega \vert \, \ll \, m $.
\end{lem}

 \begin{proof}  Note that for each
  $x^* \in \, X^*$,      the map  $\omega
\rightarrow  x^* \mu_\omega  \ \  (\Omega \rightarrow  M[0,1])$ is
weak* measurable and define an operator $ T^{x^*}: L^1 [0,1]
\rightarrow L^1 (\lambda ) $ which is nearly representable;  in fact
$T^{x^*} $ is the composition of the nearly representable operator T
with the operator $V^{x^*}:  L^1 (\lambda,X) \rightarrow L^1
(\lambda)\, (f \rightarrow x^*f ) .$ Using the fact that $L^1 (\lambda) $
has the NRNP, the operator $T^{x^*} $ is a representable
operator and therefore for
$\lambda $ a.e $\omega$, we get by Proposition~1 of \cite{FA}
that $ \vert x^*  \mu_\omega \vert \ \ll m .$
   Now using the same argument as in Lemma 2 of \cite{RS1}, we have the
conclusion of the lemma.
\end{proof}

As a consequence of Lemma~4, there exists $\Omega^\prime $ a
measurable subset of $\Omega$ with  $\lambda(\Omega \setminus
\Omega^\prime)=0$ and such that for each $\omega \in \Omega^\prime$, $\vert
\mu_\omega \vert \ll m$. Let  $g_\omega  \in  L^1 [0,1] $
 be the Radon-Nikodym density of $\vert \mu_\omega \vert $
with respect to $m $  for $\omega \in \Omega^\prime$ and $g_\omega=0 $
for  $\omega \in \Omega \setminus \Omega^\prime$.
 By $(\alpha) $
of Theorem~1,  we have the following:
for every $I$ measurable subset of $[0,1]$, the map
           $ \omega \rightarrow |\mu_\omega|(I) =\dint_{I} g_\omega(t)\ dt $
is measurable so one can deduce from the Pettis-measurability theorem
that $\omega \to g_\omega$ ($\Omega \to L^1[0,1]$) is norm-measurable.
 Moreover, $\dint_\Omega ||g_\omega||\
d\lambda(\omega) \leq ||T|| $.
  From this,   one can find a function $\Gamma \in L^1(\lambda \otimes
\mu)$  with $\Gamma(\omega,.) = g_\omega $ for $\lambda$- a.e $\omega
\in
\Omega$.

Let $V_n $ be the measurable subset of $ \Omega \times [0,1]  $ given by
   $ V_n = \left\{ (\omega,t ) ;\  n-1  \leq  \Gamma(\omega ,t)< n \right \}$.
   The $V_n $ 's are clearly disjoint and
 $\Omega\times [0,1]  =  \bigcup_n V_n$.

 Notice that for $\omega \in \Omega^\prime$,
 $ \vert \mu_\omega \vert \ll m$ and
we have $\chi_{V_n}(\omega,.) \Gamma(\omega,.)  \in L^\infty [0,1]$ and
therefore
for every $h \in L^1[0,1]$, $\chi_{V_n}(\omega,.)h(.) \Gamma(\omega ,.)
\in L^1[0,1]$ that is $\chi_{V_n}(\omega,.)h(.) \in L^1(\vert \mu_\omega
\vert)$. Hence the following map is well defined:
\begin{equation*}\begin{split}
   k_n: \Omega &\longrightarrow \cal L (L^1[0,1],X) \cr
                  \omega &\longrightarrow
\begin{cases}
k_n(\omega)(h)
=\dint \chi_{V_n}(\omega,t)h(t)  d\mu_\omega(t) \quad &\text{if
$\omega \in \Omega^\prime$}\cr
0 \quad &\text{otherwise}. \end{cases} \end{split}\end{equation*}
     It is clear that $\Vert k_n(\omega) \Vert \leq n$ for every
$\omega$.

\noindent
 Claim:  The map
$\omega
\rightarrow k_n(\omega)$ is strongly measurable:

 To prove the claim, notice that since $\sup\limits_{\omega \in
\Omega}
\Vert k_n(\omega) \Vert \leq
n$,
it is
enough to show using the denseness of the  simple functions and the
Pettis measurability theorem that for every measurable subset
$I$ of
$[0,1]$  and $x^* \in X^*$,
$\omega
\to \langle k_n(\omega)\chi_I,x^* \rangle $ is measurable;

 Let $h_\omega: [0,1] \rightarrow X^{**}$ be a
weak$^*$-density of $\mu_\omega$ with respect to $m$ for $\omega \in
\Omega^\prime $ and $0$ otherwise.  The map $\omega \rightarrow
\langle h_\omega (.), x^* \rangle $ belongs to $L^1(\lambda, L^1 [0,1])$
and  we can find a map $h \in L^1(\Omega \times [0,1])$ so that for
a.e $\omega \in \Omega$, $h(\omega,.) = h_\omega (.) $. Now the map
             $ (\omega,t)\rightarrow \chi_{V_n}(\omega,t) h(\omega,t)$
 ($\Omega \times [0,1] \to \R$)
   is  measurable and therefore for every $x^* \in X^*$
and a measurable subset $I$ of [0,1] we have
$$\omega \to \int_I \chi_{V_n}(\omega,t) \langle h_\omega (t),x^* \rangle
\ dm(t) = \langle k_n(\omega)\chi_I ,x^* \rangle. $$
This shows that $\omega \to k_n(\omega)\chi_I$ is measurable.

Let us now define  operator $ T^{(n)}: L^1[0,1] \rightarrow L^\infty(\lambda,
X)$ by $ T^{(n)}(f)= k_n(.)(f)$ and
 consider the family of measures in $ M([0,1],X)$,
$(\nu_\omega)_{\omega \in \Omega}$ defined by:
  $$ \nu_\omega (A) = \int_A \chi_{V_n}(\omega,t) \ d\mu_\omega(t) \ \
\text{for every $A$ measurable}.$$
It is clear that $k_n(\omega)(f)=\dint f(t)\ d\nu_\omega(t)$ for every $f
\in L^1[0,1]$ and
  $d|\nu_\omega| \leq n dt$.

\begin{lem} For evry $n \in \N$,
 the operator
 $J_X \circ  T^{(n)}$ is nearly representable.
 \end{lem}

\begin{proof} Let us  fix  a Dunford Pettis operator $D$ and
let $ \gamma_{k}^{(n)} = \sum\limits_{j=1}^{j_k} f_{j,k} \otimes h_{j,k} $
 be an approximating sequence for $\chi_{V_n}$ in $L^1([0,1]\times \Omega)$
with $0 \leq \gamma_{k}^{(n)} \leq \chi_{V_n}$ for every $k \in \N$
( see \cite{DS}, p.198).
Consider the sequence of operators
  $  T_{k}^{(n)}:  L^1[0,1]  \rightarrow L^1(\lambda,X) $
defined by: $$ T_{k}^{(n)} (f)(\omega) = \int \gamma_{k}^{(n)} (\omega, t)
f(t)
\ d\mu_\omega(t).$$
We claim that the operator $ T_{k}^{(n)}$ is nearly representable. Indeed,
if we
denote by $M_{f_{j,k}}$ and $M_{h_{j,k}}$ the multiplication by $f_{j,k}$
and $h_{j,k}$ respectively, we have
$ T_{k}^{(n)}
= \sum\limits_{j=1}^{j_k} M_{f_{j,k}}\circ T
\circ M_{h_{j,k}}$.
For that, let $f \in L^1[0,1]$; for a.e $\omega \in \Omega$,
 \begin{equation*}\begin{split}
\left (  \sum_{j=1}^{j_k} M_{f_{j,k}} \circ T \circ M_{h_{j,k}} \right )
(f)(\omega)  &= \sum_{j=1}^{j_k} f_{j,k} (\omega) \ T(h_{j,k}. f)(\omega)
\cr
&= \sum_{j=1}^{j_k} f_{j,k} (\omega) \ \int h_{j,k}(t) f(t)\
d\mu_\omega(t) \cr
&= \int \left ( \sum_{j=1}^{j_k}  f_{j,k}(\omega) h_{j,k}(t) f(t) \
\right  ) d\mu_\omega (t) \cr
&= \int \gamma_{k}^{(n)} (\omega,t) f(t) \ d\mu_\omega(t).
\end{split}\end{equation*}
Now since for every $j \leq j_k$,
$M_{f_{j,k}} \circ T \circ M_{h_{j,k}} \circ D$
 is representable, so is
$T_{k}^{(n)}\circ D$. To conclude the proof of the lemma, let  $\omega
\to
\nu_{k,\omega}^{D}$ and $\omega \to \nu_{\omega}^{D}$
 be the representation given by Theorem~1
 of $ T_{k}^{(n)}\circ D$ and $J_X \circ T^{(n)}\circ D$ respectively.
 We have the following:
\begin{equation*}\begin{split}
    \int |\nu_{k,\omega}^D -\nu_{\omega}^D|\ d\lambda(\omega)
     &= \int_{\Omega} \sup_{l \in \N} \sum_{m=1}^{2^l} ||\nu_{k,\omega}^D
(I_{l,m}) - \nu_{\omega}^D (I_{l,m})|| \ d\lambda(\omega) \cr
 &= \int_\Omega
\sup_{l \in \N} \sum_{m=1}^{2^l} ||\int \left (
\gamma_{k}^{(n)}(\omega,t)
-\chi_{V_n}(\omega,t)\right  )  D(\chi_{I_{l.m}})(t)  d\mu_\omega(t) ||
\ d\lambda(\omega) \cr
 &\leq \iint |\gamma_{k}^{(n)} (\omega,t) - \chi_{V_n}(\omega,t)| \
|D|(\chi_{[0,1]})(t) \ \Gamma(\omega,t)\ dt  d\lambda(\omega)
 \end{split}\end{equation*}
where $|D|$ is the modulus of $D$ (see \cite{LT}).
Notice that since $0 \leq \gamma_{k}^{(n)} \leq \chi_{V_n}$, we have
\begin{equation*}\begin{split}
 |\gamma_{k}^{(n)}(\omega,t) - \chi_{V_n}(\omega,t)|\
|D|(\chi_{[0,1]})(t) \ \Gamma(\omega,t)
 &\leq 2 \ \chi_{V_n}(\omega,t) \ |D|(\chi_{[0,1]})(t)\ \Gamma(\omega,t)
\cr
 &\leq 2n |D|(\chi_{[0,1]})(t).
\end{split}\end{equation*}
And by the Lebesgue dominated convergence,
   $\lim\limits_{k \to \infty} \dint | \nu_{k,\omega}^D -\nu_{\omega}^D| \
d\lambda(\omega) =0 $
and hence by passing to a subsequence (if necessary), we may assume that
  $\lim\limits_{k \to \infty}|\nu_{k,\omega}^D
-\nu_{\omega}^D|=0$ for a.e $\omega \in \Omega$.

Fix $B_0$ a subset of $\Omega$ with $\lambda(B_0) =0 $ and for every
$\omega \notin B_0$, $\lim\limits_{k \to \infty} |\nu_{k,\omega}^D -
\nu_{\omega}^D | =0$. Since $ T_{k}^{(n)} \circ D$ is
representable,
one can find a subset $B_k$ of $\Omega$ with $\lambda(B_k)=0$ and such
that for each $\omega \notin B_k$, $\nu_{k,\omega}^D$ has Bochner
integrable density. We can conclude that  for $\omega \notin
\bigcup_{k=0}^\infty B_k $, the measure $\nu_{\omega}^D$ is the limit for
the variation norm of a sequence of measures with Bochner integrable
densities and therefore has Bochner integrable
density. Now using Proposition~1, the operator
$J_X \circ T^{(n)}\circ D$ is representable. The lemma is proved.
\end{proof}

We are now ready to  complete the
 proof of the theorem:
By Proposition~5, the operator $J_X\circ T^{(n)}$ is representable
and therefore the operator $K_n: L^1(\Omega \times [0,1]) \rightarrow
X$ given by $K_n(f)= \dint k_n(\omega)(f(\omega,.))\ d\lambda(\omega)$
is representable (see Proposition~4).

Let $\phi_n: \Omega \times [0,1] \rightarrow X$ be a representation
of
 $K_n$ and consider
 $     \varphi=\sum\limits_{n=1}^ \infty \phi_n \ \chi_{V_n} $.

\noindent
We claim that $\varphi$ belongs to $ L^1(\Omega \times [0,1],X) $.

For that, fix $\alpha_\omega: [0,1] \rightarrow X^{**}$ a weak$^*$-
density of $\mu_\omega$ with respect to $\vert \mu_\omega \vert$(see
\cite{DIU} or \cite{IT}). Clearly,
$\Vert \alpha_\omega (t) \Vert =1$ for a.e $t \in [0,1]$ and
$d\mu_\omega= \alpha_\omega \Gamma(\omega,.) dt$.
Let $G_n: \Sigma_{\lambda \otimes m} \rightarrow X$ be given by
 $G_n(V)= K_n(\chi_V)$.

\noindent By the definition of $K_n$, $G_n(V) =
\text{weak}^*-\displaystyle\iint_{V}
\chi_{V_n}(\omega,t) \ \alpha_{\omega}(t)\ \Gamma(\omega,t) \ dt\
d\lambda(\omega)$. In the other hand since $K_n$ is represented by
$\phi_n$, we have $G_n(V)= \displaystyle\iint_{V} \phi_n(\omega,t) \ dt \
d\lambda(\omega)$. So we have
 $$ || \phi_n || =
\vert G_n \vert (\Omega \times [0,1]) =\iint\Vert\phi_n(\omega,t)
\Vert \ d\lambda\otimes m (\omega,t)$$
and using the weak$^*$-density, we get
$$ \Vert \phi_n \Vert =\iint\chi_{V_n}(\omega,t)
\Gamma(\omega,t) d\lambda\otimes m(\omega,t)$$
which shows that  $\sum\limits_{n=1}^\infty || \phi_n \chi_{V_n} ||_1  \leq
 \displaystyle\iint \Gamma(\omega,t) \ dt \ d\lambda(\omega)$.
Hence the series is  convergent.

For each $V \in \Sigma_{\lambda \otimes m} $, we get
\begin{equation*}\begin{split}
  \iint_{V} d\mu_\omega (t) \ d\lambda (\omega)
    &=\sum_{n=1}^\infty \iint_V \chi_{V_n}(\omega,t) \ d\mu_\omega (t) \
d\lambda(\omega) \cr
    &=\sum_{n=1}^\infty K_n(\chi_V) = \sum_{n=1}^\infty K_n (\chi_V .
\chi_{V_n}) \cr
    &=\sum_{n=1}^\infty \iint_V \phi_n (\omega,t) \ dt\
d\lambda(\omega) \cr
   &=\iint_{V} \varphi (\omega,t) \ dt\ d\lambda(\omega).
   \end{split}\end{equation*}
In particular, for each $ \ A \in \Sigma_m $  and $ \  B \in
\Sigma_\lambda$,
 $\dint_{B} \mu_\omega(A)\ d\lambda (\omega)
=\dint_{B} \{ \dint_{A} \varphi (\omega,t) \ dt \}\ d\lambda
(\omega) $
which shows that  $\mu_\omega (A) =\dint_{A} \varphi (t,\omega)\ dt$
 for a.e  $\omega$. The theorem is proved . $\qed$

\medskip
Before stating the next  extension, let us recall  (as in
\cite{T2}) that, if E is a K\"othe function space on
 $(\Omega, \Sigma, \lambda)$ (in the sense of \cite{LT})
and $X$ is a Banach space then   $E(X)$ will  be the space of
all (classes of) measurable map $f: \Omega \rightarrow X $  so that
$\omega \rightarrow \Vert f(\omega) \Vert $  belongs to $E$.

\noindent
{\bf Corollary.}  {\it If E does not contain a copy of  $c_0 $  and X
has the NRNP, then E(X) has the NRNP.}

\begin{proof}    Without loss of generality, we can assume that
$E$ is  order continuous,  $(\Omega, \Sigma, \lambda)$
is a separable probability space (see \cite{LT}) and the Banach space $X$
is separable. By a result of Lotz, Peck and Porta (\cite{LPP}),
the inclusion map from $E$ into $L^1(\lambda)$ is a semi-imbedding.
The same is true for the inclusion $J_X: E(X) \to L^1(\lambda,X)$ (see
 \cite{RS3} Lemma~3).
Now let  $T: L^1 [0,1] \rightarrow E(X) $ be a nearly representable
operator.
The operator $ J_X \circ T $ is also nearly representable and hence
representable
(by Theorem~2). So the operator T must be representable (see \cite{BR}).
\end{proof}

\section{Concluding remarks}
If $X$ and $Y$ are Banach spaces with the NRNP, then
 $X \widehat{\otimes}_\pi Y$ ( $\widehat{\otimes}_\pi$ is the projective tensor
product) need not satisfy the NRNP. This can be seen from Pisier's famous
example that $L^1/{H_{0}^1} \widehat{\otimes}_\pi L^1/{H_{0}^1}$ contains $c_0$
(hence failing  the NRNP) while $L^1/{H_{0}^1}$ has the NRNP.

If $X$ is a Banach space and $(\Omega, \Sigma)$ is a measure space, we denote
by $M(\Omega,X^*)$ the space of $X^*$-valued $\sigma$-additive measures 
of bounded variation with the usual total variation norm. In light of
Theorem~2, one can ask the following question: Does $M(\Omega,X^*)$ have the 
NRNP whenever $X^*$ does ? 
 It should be noted that for non-dual space, the answer is negative: the space
$E$ constructed by Talagand in \cite{T5} is a Banach lattice that does not
contain $c_0$ (so it has the NRNP) but $M(\Omega, E)$ contains $c_0$.

Finally, since $L^1$-spaces are the primary examples of Banach spaces with 
the NRNP,  the following question arises: Do non-commutative $L^1$-spaces
have the NRNP?
Note that since $C_1$ (the trace class operators) has the RNP, it has the
NRNP; however it is still unknown if $C_E$ has the NRNP if $E$ is a 
symmetric sequence space that does not contain $c_0$. We remark that
non-commutative $L^1$-spaces have the ARNP (\cite{HP}).
   

\end{document}